\documentclass[journal]{IEEEtran}%
\usepackage{amssymb}
\usepackage{amsmath}
\usepackage{amsfonts}%
\usepackage{graphicx}
\providecommand{\U}[1]{\protect\rule{.1in}{.1in}}
\newtheorem{theorem}{Theorem}

\newtheorem{definition}[theorem]{Definition}
\newtheorem{example}[theorem]{Example}

\newtheorem{lemma}[theorem]{Lemma}

\newtheorem{proposition}[theorem]{Proposition}
\newtheorem{remark}[theorem]{Remark}

\begin{document}

\title{Motion Planning for Kinematic systems}
\author{Nicolas Boizot and Jean-Paul Gauthier \thanks{The authors are with LSIS,
laboratoire des sciences de l'information et des syst\`{e}mes, UMR\ CNRS 6168,
Domaine universitaire de Saint J\'{e}r\^{o}me, Avenue Escadrille Normandie
Niemen, 13397 MARSEILLE Cedex 2, France. \ J.P.\ Gauthier is also with INRIA
team GECO.}}
\pubid{0000--0000/00\$00.00~\copyright ~2012 IEEE}
\maketitle

\begin{abstract}
\textrm{In this paper, we present a general theory of motion planning for
kinematic systems.\ This theory has been developed for long by one of the
authors in a previous series of papers.\ It is mostly based upon concepts from
subriemannian geometry.\ Here, we summarize the results of the theory, and we
improve on, by developping in details an intricated case: the ball with a
trailer, which corresponds to a distribution with flag of type 2,3,5,6.}

\textrm{\ This paper is dedicated to Bernard Bonnard for his 60}$^{th}$ birthday.

\end{abstract}


\tableofcontents

\section{Introduction}

Here we present the main lines of a theory of motion planning for kinematic
systems, which is developped for about ten years in the papers
\ \cite{GM,GZ,GZ2,WIL,GZ3,ano,jak}.\ One of the purposes of the paper is to
survey the whole theory disseminated in these papers. But also we improve on
the theory, by treating one more case, in which "the fourth order brackets are
involved". We aslo improve on several previous results (periodicity of our
optimal trajectories for instance).\ Potential application of this theory is
motion planning for kinematic robots.\ We will show several basic examples here.

The theory starts from the seminal work of F.\ Jean, in the papers
\cite{J1,J2,J3}.\ At the root of this point of view in robotics, there are
also more applied authors like J.P.\ Laumond \cite{lau}. See also \cite{liu}.

We consider kinematic systems that are given under the guise of a
vector-distribution $\Delta$ over a $n$-dimensional manifold $M$. The rank of
the distribution is $p$, and the corank $k=n-p.\ $Motion planning problems
will aways be local problems in an open neighborhood of a given finite path
$\Gamma$ in $M.$ Then we may always consider that $M=\mathbb{R}^{n}.$ From a
control point of view, a kinematic system can be specified by a control
system, linear in the controls, typically denoted by $\Sigma$:%
\begin{equation}
(\Sigma)\text{ }\dot{x}=%
{\displaystyle\sum\limits_{i=1}^{p}}
F_{i}(x)u_{i},\label{sys1}%
\end{equation}

where the $F_{i}$'s are smooth ($C^{\infty})$ vector fields that span the
distribution $\Delta.$ The standard controllability assumption is always
assumed, i.e.\ the Lie algebra generated by the $F_{i}$'s is transitive on
$M.$ Consequently, the distribution $\Delta$ is \textit{completely
nonintegrable}, and any smooth path $\Gamma:[0,T]\rightarrow M$ can be
unifomly approximated by an admissible path $\gamma:[0,\theta]\rightarrow M $,
i.e. a Lipschitz path, which is almost everywhere tangent to $\Delta,$ i.e., a
trajectory of (\ref{sys1})$.$

This is precisely the \textit{abstract answer} to the kinematic motion
planning probem: \textit{it is possible to approximate uniformly nonadmissible
paths by admissible ones}. The purpose of this paper is to present a general
constructive theory that solves this problem in a certain \textit{optimal} way.

\bigskip More precisely, in this class of problems, it is natural to try to
minimize a cost of the following form:%
\[
J(u)=%
{\displaystyle\int\limits_{0^{{}}}^{\theta}}
\sqrt{%
{\displaystyle\sum\limits_{i=1}^{p}}
(u_{i})^{2}}dt,
\]
for several reasons: 1.\ the optimal curves do not depend on their
parametrization, 2.\ the minimization of such a cost produces a metric space
(the associated distance is called the subriemannian distance, or the
Carnot-Caratheodory distance), 3.\ Minimizing such a cost is equivalent to
minimize the following (called the \textit{energy} of the path) quadratic cost
$J_{E}(u)$, in fixed time $\theta$:%
\[
J_{E}(u)=%
{\displaystyle\int\limits_{0^{{}}}^{\theta}}
{\displaystyle\sum\limits_{i=1}^{p}}
(u_{i})^{2}dt.
\]

The distance is defined as the minimum length of admissible curves connecting
two points, and the length of the admissible curve corresponding to the
control $u:[0,\theta]\rightarrow M$ is just $J(u).$

In this presentation, another way to interpret the problem is as follows: the
dynamics is specified by the distribution $\Delta$ (i.e. not by the vector
fields $F_{i},$ but their span only).\ The cost is then determined by an
Euclidean metric $g$ over $\Delta,$ specified here by the fact that the
$F_{i}$'s \ form an orthonormal frame field for the metric.

At this point we would like to make a more or less philosophical comment:
there is, in the world of nonlinear control theory, a permanent twofold critic
against the optimal control approach: 1. the choice of the cost to be
minimized is in general rather arbitrary, and 2.\ optimal control solutions
may be non robust.\ 

Some remarkable conclusions of our theory show the following: in reasonable
dimensions and codimensions, the optimal trajectories are extremely robust,
and in particular, do not depend at all (modulo certain natural
transformations) on the choice of the metric, but on the distribution $\Delta$
only. Even stronger: they depend only on the \textit{nilpotent approximation
along} $\Gamma$ (a concept that will be defined later on, which is a good
local approximation of the problem).\ For a lot of low values of the rank $p$
and corank $k,$ these nilpotent approximations have no parameter (hence they
are in a sense universal).\ The \textit{asymptotic optimal sysntheses} (i.e.
the phase portraits of the admissible trajectories that approximate up to a
small $\varepsilon)$ are also universal.

\bigskip

Given a motion planning problem, specified by a (nonadmissible) curve
$\Gamma,$ and a Subriemannian structure (\ref{sys1}), we will consider two
distinct concepts, namely: 1.\ The \textit{metric complexity} $MC(\varepsilon
)$ that measures asymptotically the length of the best $\varepsilon
$-approximating admissible trajectories$,$ and 2.\ The \textit{interpolation
entropy} $E(\varepsilon)$, that measures the length of the best admissible
curves that interpolate $\Gamma$ with pieces of length $\varepsilon.$

The first concept was introduced by F.\ Jean in his basic paper \cite{J1}. The
second concept is closely related with the entropy of F.\ Jean in \cite{J2},
which is more or less the same as the Kolmogorov's entropy of the path
$\Gamma,$ for the metric structure induced by the Carnot-Caratheodory metric
of the ambient space.

\bigskip Also, along the paper, we will deal with \textit{generic} problems
only (but generic in the global sense, i.e.\ stable singularities are
considered). That is, the set of motion planning problems on $\mathbb{R}^{n}$
is the set of couples $(\Gamma,\Sigma),$ embedded with the $C^{\infty}$
topology of uniform convergence over compact sets, and generic problems (or
\textit{problems in general position}) form an open-dense set in this
topology. For instance, it means that the curve $\Gamma$ is always tranversal
to $\Delta$ (except maybe at isolated points, in the cases $k=1$ only).
Another example is the case of a surface of degeneracy of the Lie bracket
distribution $[\Delta,\Delta]$ in the $n=3,$ $k=1$ case. Generically, this
surface (the Martinet surface) is smooth, and $\Gamma$ intersects it
transversally at a finite number of points only.

Also, along the paper, we will illustrate our results with one of the
following well known academic examples:

\begin{example}
\label{unic}the unicycle:
\begin{equation}
\dot{x}=\cos(\theta)u_{1},\text{ }\dot{y}=\sin(\theta)u_{1},\text{ }%
\dot{\theta}=u_{2}\label{unicycle}%
\end{equation}

\end{example}

\begin{example}
\label{ctrl}\bigskip the car with a trailer:%
\begin{equation}
\dot{x}=\cos(\theta)u_{1},\text{ }\dot{y}=\sin(\theta)u_{1},\text{ }%
\dot{\theta}=u_{2},\text{ }\dot{\varphi}=u_{1}-\sin(\varphi)u_{2}%
\label{cartrailer}%
\end{equation}

\end{example}

\begin{example}
\label{bpln}\bigskip the ball rolling on a plane:%
\begin{equation}
\dot{x}=u_{1},\text{ }\dot{y}=u2,\text{ }\dot{R}=\left[
\begin{array}
[c]{ccc}%
0 & 0 & u_{1}\\
0 & 0 & u_{2}\\
-u_{1} & -u_{2} & 0
\end{array}
\right]  R,\label{brp}%
\end{equation}
where $(x,y)$ are the coordinates of the contact point between the ball and
the plane, $R\in SO(3,\mathbb{R})$ is the right orthogonal matrix representing
an othonormal frame attached to the ball.
\end{example}

\begin{example}
\label{brpt}the ball with a trailer%
\begin{align}
\dot{x}  & =u_{1},\text{ }\dot{y}=u2,\text{ }\dot{R}=\left[
\begin{array}
[c]{ccc}%
0 & 0 & u_{1}\\
0 & 0 & u_{2}\\
-u_{1} & -u_{2} & 0
\end{array}
\right]  R,\label{balltrailer}\\
\text{ }\dot{\theta}  & =-\frac{1}{L}(\cos(\theta)u_{1}+\sin(\theta
)u_{2}).\nonumber
\end{align}

\end{example}

\begin{figure}[ptb]%
\centering
\includegraphics[
natheight=1.802300in,
natwidth=1.666500in,
height=1.8395in,
width=1.7028in
]%
{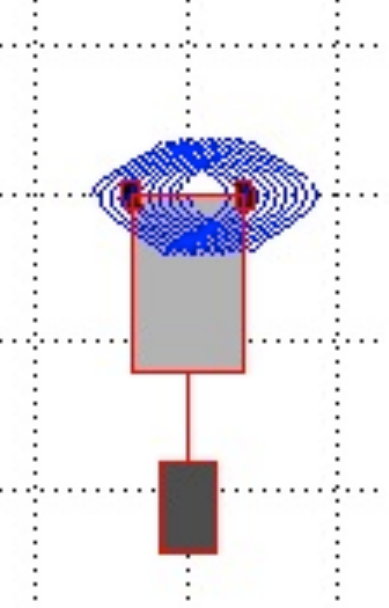}%
\caption{Parking of the car with a trailer}%
\label{fig1}%
\end{figure}

\begin{figure}[ptb]%
\centering
\includegraphics[
natheight=2.739700in,
natwidth=3.312200in,
height=2.6965in,
width=3.2543in
]%
{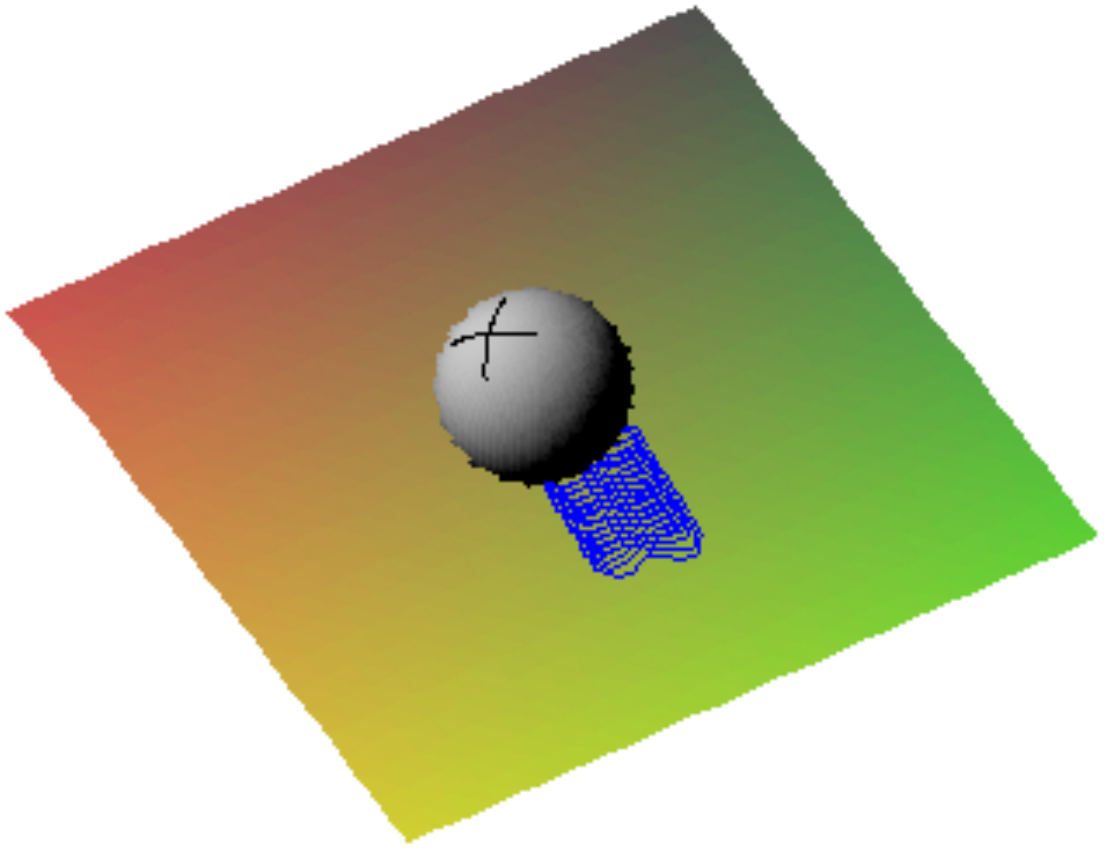}%
\caption{Approximating rolling with slipping}%
\label{fig2}%
\end{figure}

Typical motion planning problems are: 1.\ for example (\ref{ctrl}%
),\textit{\ the parking problem}: the non admissible curve $\Gamma$ is
$s\rightarrow(x(s),y(s),\theta(s),\varphi(s))=(s,0,\frac{\pi}{2},0),$ 2.\ for
example (\ref{bpln}), the \textit{full rolling with slipping problem},
$\Gamma:s\rightarrow(x(s),y(s),R(s))$ $=(s,0,Id),$ where $Id$ is the identity
matrix. On figures \ref{fig1}, \ref{fig2} we show our approximating
trajectories for both problems, that are in a sense universal.\ In figure
\ref{fig1}, of course, the $x$-scale is much larger than the $y$-scale.

\bigskip

\bigskip%

%

Up to now, our theory covers the following cases:

(C1) The distribution $\Delta$ is one-step bracket generating (i.e.
$\dim([\Delta,\Delta]=n)$ except maybe at generic singularities,

(C2) The number of controls (the dimension of $\Delta)$ is $p=2,$ and
$n\leq6.$

The paper is organized as follows: In the next section \ref{prereq}, we
introduce the basic concepts, namely the metric complexity, the interpolation
entropy, the nilpotent approximation along $\Gamma,$ and the normal
coordinates, that will be our basic tools.

Section \ref{results} summarizes the main results of our theory, disseminated
in our previous papers, with some complements and details. Section \ref{new}
is the detailed study of the case $n=6,$ $k=4,$ which corresponds in
particular to example \ref{brpt}, the ball with a trailer. In Section
\ref{concl}, we state a certain number of remarks, expectations and conclusions.

\section{Basic concepts \label{prereq}}

In this section, we fix a generic motion planning problem $\mathcal{P=}%
(\Gamma,\Sigma).$ Also, along the paper there is a small parameter
$\varepsilon$ (we want to approximate up to $\varepsilon),$ and certain
quantities $f(\varepsilon),g(\varepsilon)$ go to $+\infty$ when $\varepsilon$
tends to zero.\ We say that such quantities are equivalent $(f\simeq g)$ if
$\lim_{\varepsilon\rightarrow0}\frac{f(\varepsilon)}{g(\varepsilon)}=1.$ Also,
$d$ denotes the subriemannian distance, and we consider the $\varepsilon
$-subriemammian tube $T\varepsilon$ and cylinder $C\varepsilon$ around
$\Gamma:$%

\begin{align*}
T_{\varepsilon}  & =\{x\in M\text{ }|\text{ }d(x,\Gamma)\leq\varepsilon\},\\
C_{\varepsilon}  & =\{x\in M\text{ }|\text{ }d(x,\Gamma)=\varepsilon\}.
\end{align*}

\subsection{Entropy versus metric complexity \label{entcomp}}

\begin{definition}
\label{mc}The \textit{metric complexity} $MC(\varepsilon)$ of $\mathcal{P}$ is
$\frac{1}{\varepsilon}$ times the minimum length of an admissible curve
$\gamma_{\varepsilon}$ connecting the endpoints $\Gamma(0),$ $\Gamma(T)$ of
$\Gamma,$ and remaining in the tube $T_{\varepsilon}.$
\end{definition}

\begin{definition}
\label{mce}The \textit{interpolation entropy} $E(\varepsilon)$ of
$\mathcal{P}$\ is $\frac{1}{\varepsilon}$ times the minimum length of an
admissible curve $\gamma_{\varepsilon}$ connecting the endpoints
$\Gamma(0),\Gamma(T)$ of $\Gamma,$ and \ $\varepsilon$-interpolating $\Gamma$,
that is, in any segment of $\gamma_{\varepsilon}$ of length $\geq\varepsilon,$
there is a point of $\Gamma.$
\end{definition}

These quantities $MC(\varepsilon),E(\varepsilon)$ are functions of
$\varepsilon$ which tends to $+\infty$ as $\varepsilon$ tends to zero. They
are considered \textbf{up to equivalence}.The reason to divide by
$\varepsilon$ is that the second quantity counts the number of $\varepsilon
$-balls to cover $\Gamma,$ or the number of pieces of length $\varepsilon$ to
interpolate the full path. This is also the reason for the name "entropy".

\begin{definition}
An asymptotic optimal synthesis is a one-parameter family $\gamma
_{\varepsilon}$ of admissible curves, that realizes the metric complexity or
the entropy.
\end{definition}

Our main purpose in the paper is twofold:

1.\ We want to estimate the metric complexity and the entropy, in terms of
certain invariants of the problem.\ Actually, in all the cases treated in this
paper, we will give eplicit formulas.

2.\ We shall exhibit explicit asymptotic optimal syntheses realizing the
metric complexity or/and the entropy.

\subsection{Normal coordinates\label{normc}}

\ Take a \textbf{parametrized} $p$-dimensional surface $S,$ transversal to
$\Delta$ (maybe defined in a neighborhood of $\Gamma$ only)$,$
\[
S=\{q(s_{1},...,s_{p-1},t)\in\mathbb{R}^{n}\},\text{with }q(0,...,0,t)=\Gamma
(t).
\]
Such a \textit{germ} exists if $\Gamma$ is not tangent to $\Delta.$ The
exclusion of a neighborhood of an isolated point where $\Gamma$ is tangent to
$\Delta$, (that is $\Gamma$ becomes \textquotedblright almost
admissible\textquotedblright), will not affect our estimates presented later
on (it will provide a term of higher order in $\varepsilon).$\ \ .

\smallskip In the following, $\mathcal{C}_{\varepsilon}^{S}$ will denote the
cylinder $\{\xi;$ $d(S,\xi)=\varepsilon\}.$

\begin{lemma}
\label{nco} (Normal coordinates with respect to $S).$ There are mappings
$x:\mathbb{R}^{n}\rightarrow\mathbb{R}^{p},$ $y:\mathbb{R}^{n}\rightarrow
\mathbb{R}^{k-1},$ $w:\mathbb{R}^{n}\rightarrow\mathbb{R},$ such that
$\xi=(x,y,w)$ is a coordinate system on some neighborhood of $S$ in
$\mathbb{R}^{n}$, such that:

0. $S(y,w)=(0,y,w),$ $\Gamma=\{(0,0,w)\}$

1. The restriction\ $\Delta_{|S}=\ker dw\cap_{i=1,..k-1}\ker dy_{i},$ the
metric $g_{|S}=\sum_{i=1}^{p}(dx_{i})^{2},$

2. $\mathcal{C}_{\varepsilon}^{S}=\{\xi|\sum_{i=1}^{p}x_{i}{}^{2}%
=\varepsilon^{2}\},$

3. geodesics of the Pontryagin's maximum principle (\cite{PMP}) meeting the
transversality conditions w.r.t. $S$ are the straight lines through $S,$
contained in the planes $P_{y_{0},w_{0}}=\{\xi|(y,w)=(y_{0},w_{0})\}.$ Hence,
they are orthogonal to $S.$

These normal coordinates are unique up to changes of coordinates of the form
\begin{equation}
\tilde{x}=T(y,w)x,(\tilde{y},\tilde{w})=(y,w),\label{ccor}%
\end{equation}
where $T(y,w)\in O(p),$ the $p$-orthogonal group.
\end{lemma}

\subsection{Normal forms, Nilpotent approximation along $\Gamma$\label{nform}}

\subsubsection{Frames\label{fra}}

Let us denote by $F=(F_{1},...,F_{p})$ the orthonormal frame of vector fields
generating $\Delta.$ Hence, we will also write $\mathcal{P}=(\Gamma,F).$ If a
global coordinate system $(x,y,w)$, not necessarily normal, is given on a
neighborhood of $\Gamma$ in $\mathbb{R}^{n},$ with $x\in\mathbb{R}^{p},$
$y\in\mathbb{R}^{k-1},$ $w\in\mathbb{R},$ then we write:%
\begin{align}
F_{j}  & =\sum_{i=1}^{p}\mathcal{Q}_{i,j}(x,y,w)\frac{\partial}{\partial
x_{i}}+\sum_{i=1}^{k-1}\mathcal{L}_{i,j}(x,y,w)\frac{\partial}{\partial y_{i}%
}\label{QLM}\\
& +\mathcal{M}_{j}(x,y,w)\frac{\partial}{\partial w},\text{ }\nonumber\\
\text{\ \ \ }j  & =1,...,p.\nonumber
\end{align}
Hence, the SR metric is specified by the triple $(\mathcal{Q},\mathcal{L}%
,\mathcal{M})$ of smooth $x,y,w$-dependent matrices.

\subsubsection{The general normal form\label{gennf}}

\smallskip Fix a surface $S$ as in Section \ref{ncor} and a normal coordinate
system $\xi=(x,y,w)$ for a problem $\mathcal{P}.$

\begin{theorem}
\label{normal}(Normal form, \cite{AG2}) There is a unique orthonormal frame
$F=(\mathcal{Q},\mathcal{L},\mathcal{M})$ for ($\Delta,g)$ with the following properties:

1.\ $\mathcal{Q}(x,y,w)$ is symmetric, $\mathcal{Q}(0,y,w)=Id$ (the identity matrix),

2. $\mathcal{Q}(x,y,w)x=x,$

3. $\mathcal{L}(x,y,w)x=0,$ $\mathcal{M}(x,y,w)x=0.$

4. Conversely if $\xi=(x,y,w)$ is a coordinate system satisfying conditions 1,
2, 3 above, then $\xi$ is a normal coordinate system for the SR metric defined
by the orthonormal frame $F$ with respect to the parametrized surface
$\{(0,y,w)\}.$
\end{theorem}

Clearly, this normal form is invariant under the changes of normal coordinates
(\ref{ccor}).

Let us write:
\begin{align*}
\mathcal{Q}(x,y,w)  &  =Id+Q_{1}(x,y,w)+Q_{2}(x,y,w)+...,\\
\mathcal{L}(x,y,w)  &  =0+L_{1}(x,y,w)+L_{2}(x,y,w)+...,\\
\mathcal{M}(x,y,w)  &  =0+M_{1}(x,y,w)+M_{2}(x,y,w)+...,
\end{align*}

where $Q_{r},L_{r},M_{r}$ are matrices depending on $\xi=(x,y,w),$ the
coefficients of which have order $k$ w.r.t. $x$ (i.e. they are in the $r^{th}
$ power of the ideal of $C^{\infty}(x,y,w)$ generated by the functions
$x_{r},$ $r=1,...,n-p).$ In particular, $Q_{1}$ is linear in $x,$ $Q_{2}$ is
quadratic, etc...\ Set $u=(u_{1},...,u_{p})\in\mathbb{R}^{p}.$ Then
$\sum_{j=1}^{k-1}L_{1_{j}}(x,y,w)u_{j}$ $=L_{1,y,w}(x,u)$ is quadratic in
$(x,u),$ and $\mathbb{R}^{k-1}$-valued. Its $i^{th}$ component is the
quadratic expression denoted by $L_{1,i,y,w}(x,u)$.\ Similarly $\sum
_{j=1}^{k-1}M_{1_{j}}(x,y,w)u_{j}$ $=M_{1,y,w}(x,u)$ is a quadratic form in
$(x,u).$ The corresponding matrices are denoted by $L_{1,i,y,w},$
$i=1,...,k-1,$ and $M_{1,y,w}.$

The following was proved in \cite{AG2}, \cite{char} for corank 1:

\begin{proposition}
\label{norprop} 1.\ $Q_{1}=0,$

2.\ $L_{1,i,y,w},$ $i=1,...,p-1,$ and $M_{1,y,w}$ are skew symmetric matrices.\ 
\end{proposition}

A first useful very rough estimate in normal coordinates is the following:

\begin{proposition}
\label{propb}\ If $\xi=(x,y,w)\in T_{\varepsilon},$ then:
\begin{align*}
||x||_{2} &  \leq\varepsilon,\\
||y||_{2} &  \leq k\varepsilon^{2},
\end{align*}
for some $k>0.$
\end{proposition}

At this point, we shall split the problems under consideration into two
distinct cases: first the 2-step bracket-generating case, and second, the
2-control case.

\subsubsection{Two-step bracket-generating case\label{spe}}

In that case, we set, in accordance to Proposition \ref{propb}, that $x$ has
weight 1, and the $y_{i}$'s and $w$ have weight 2$.$ Then, the vector fields
$\frac{\partial}{\partial x_{i}}$\smallskip\ have weight -1, and
$\frac{\partial}{\partial y_{i}},\frac{\partial}{\partial w}$ have weight
$-2.$

Inside a tube $T_{\varepsilon},$ we write our control system as a term of
order -1, plus a residue, that has a certain order w.r.t. $\varepsilon
.\ $Here, $O(\varepsilon^{k})$ means a smooth term bounded by $c\varepsilon
^{k}.$ We have, for a trajectory remaining inside $T_{\varepsilon}$:%

\begin{align}
\dot{x}  & =u+O(\varepsilon^{2});\text{ \ \ \ \ \ \ \ \ \ \ \ \ \ (1)}%
\label{estco1}\\
\dot{y}_{i}  & =\frac{1}{2}x^{\prime}L^{i}(w)u+O(\varepsilon^{2});\text{
\ \ }i=1,...,k-1;\nonumber\\
\dot{w}  & =\frac{1}{2}x^{\prime}M(w)u+O(\varepsilon^{2}),\nonumber
\end{align}

\bigskip where $L^{i}(w),M(w)$ are skew-symmetric matrices depending smoothly
on $w.\ $

\begin{remark}
In \ref{estco1}, (1), the term $O(\varepsilon^{2})$ can seem surprising.\ One
should wait for $O(\varepsilon).\ $It is due to (1) in Proposition
\ref{norprop}.
\end{remark}

In that case, we define the \textbf{Nilpotent Approximation }$\hat{P}%
$\textbf{\ along }$\Gamma$ of the problem $\mathcal{P}$ by keeping only the
term of order -1:%

\begin{align}
\dot{x}  & =u;\label{nilap1}\\
(\mathcal{\hat{P})}\text{ \ \ \ \ \ \ \ \ }\dot{y}_{i}  & =\frac{1}%
{2}x^{\prime}L^{i}(w)u;\text{ \ \ }i=1,...,p-1;\nonumber\\
\dot{w}  & =\frac{1}{2}x^{\prime}M(w)u.\nonumber
\end{align}

Consider two trajectories $\xi(t),\hat{\xi}(t)$ of $\mathcal{P}$ and
$\mathcal{\hat{P}}$ corresponding to the same control $u(t),$ issued from the
same point on $\Gamma,$ and both arclength-parametrized (which is equivalent
to $||u(t)||=1).$ For $t\leq\varepsilon,$ we have the following estimates:

\bigskip%
\begin{equation}
||x(t)-\hat{x}(t)||\leq c\varepsilon^{3},||y(t)-\hat{y}(t)||\leq
c\varepsilon^{3},||w(t)-\hat{w}(t)||\leq c\varepsilon^{3},\label{ff0}%
\end{equation}
for a suitable constant $c.\ $

\begin{remark}
\label{dist1}It follows that the distance (either $d$ or $\hat{d}$-the
distance associated with the nilpotent approximation$)$ between $\xi
(t),\hat{\xi}(t)$ is smaller than $\varepsilon^{1+\alpha}$ for some $\alpha>0.
$
\end{remark}

This fact comes from the estimate just given, and the standard ball-box
Theorem (\cite{GRO}). It will be the key point to reduce the motion planning
problem to the one of its nilpotent approximation along $\Gamma$.

\subsubsection{\bigskip The 2-control case\label{2control}}

\subsubsection{Normal forms\label{nffs}}

In that case, we have the following general normal form, in normal
coordinates. It was proven first in \cite{AmPetr92}, in the corank1 case.\ The
proof holds in any corank, without modification.

Consider Normal coordinates with respect to any surface $\mathcal{S}$. There
are smooth functions, $\beta(x,y,w),\gamma_{i}(x,y,w),\delta(x,y,w),$ such
that $\mathcal{P}$ can be written as (on a neighborhood of $\Gamma)$:

\bigskip%
\begin{align}
\dot{x}_{1}  & =(1+(x_{2}^{{}})^{2}\beta)u_{1}-x_{1}x_{2}\beta u_{2},\text{
\ }\label{nf2}\\
\text{\ }\dot{x}_{2}  & =(1+(x_{1}^{{}})^{2}\beta)u_{2}-x_{1}x_{2}\beta
u_{1},\nonumber\\
\dot{y}_{i}  & =\gamma_{i}(\frac{x_{2}}{2}u_{1}-\frac{x_{1}}{2}u_{2}%
),\text{\ \ }\dot{w}=\delta(\frac{x_{2}}{2}u_{1}-\frac{x_{1}}{2}%
u_{2}),\nonumber
\end{align}

where moreover $\beta$ vanishes on the surface $\mathcal{S}$.

The following normal forms can be obtained, on the tube $T_{\varepsilon},$ by
just changing coordinates in $\mathcal{S}$ in certain appropriate way. It
means that a trajectory $\xi(t)$ of $\mathcal{P}$ remaining in $T_{\varepsilon
}$ satisfies:

\textbf{Generic }$4-2$\textbf{\ case (see \cite{GZ3})}$:$%
\begin{align*}
\dot{x}_{1}  & =u_{1}+0(\varepsilon^{3}),\dot{x}_{2}=u_{2}+0(\varepsilon
^{3}),\\
\dot{y}  & =(\frac{x_{2}}{2}u_{1}-\frac{x_{1}}{2}u_{2})+O(\varepsilon^{2}),\\
\dot{w}  & =\delta(w)x_{1}(\frac{x_{2}}{2}u_{1}-\frac{x_{1}}{2}u_{2}%
)+O(\varepsilon^{3}).
\end{align*}
We define the nilpotent approximation as:%
\begin{align*}
(\mathcal{\hat{P}}_{4,2})\text{ \ \ }\dot{x}_{1}  & =u_{1},\dot{x}_{2}%
=u_{2},\dot{y}=(\frac{x_{2}}{2}u_{1}-\frac{x_{1}}{2}u_{2}),\\
\dot{w}  & =\delta(w)x_{1}(\frac{x_{2}}{2}u_{1}-\frac{x_{1}}{2}u_{2}%
).\text{\ }%
\end{align*}

Again, we consider two trajectories $\xi(t),\hat{\xi}(t)$ of $\mathcal{P}$ and
$\mathcal{\hat{P}}$ corresponding to the same control $u(t),$ issued from the
same point on $\Gamma,$ and both arclength-parametrized (which is equivalent
to $||u(t)||=1).$ For $t\leq\varepsilon,$ we have the following estimates:

\bigskip%
\begin{equation}
||x(t)-\hat{x}(t)||\leq c\varepsilon^{4},||y(t)-\hat{y}(t)||\leq
c\varepsilon^{3},||w(t)-\hat{w}(t)||\leq c\varepsilon^{4}.\label{ff1}%
\end{equation}

Which implies that, for $t\leq\varepsilon,$ the distance ($d$ or $\hat{d})$
between $\xi(t)$ and $\hat{\xi}(t)$ is less than $\varepsilon^{1+\alpha}$ for
some $\alpha>0,$ and this will be also the keypoint to reduce our problem to
the Nilpotent approximation.

\textbf{Generic }$5-2$\textbf{\ case (see \cite{ano})}$:$%
\begin{align*}
\dot{x}_{1}  & =u_{1}+0(\varepsilon^{3}),\dot{x}_{2}=u_{2}+0(\varepsilon
^{3}),\\
\dot{y}  & =(\frac{x_{2}}{2}u_{1}-\frac{x_{1}}{2}u_{2})+O(\varepsilon^{2}),\\
\dot{z}  & =x_{2}(\frac{x_{2}}{2}u_{1}-\frac{x_{1}}{2}u_{2})+O(\varepsilon
^{3}),\\
\dot{w}  & =\delta(w)x_{1}(\frac{x_{2}}{2}u_{1}-\frac{x_{1}}{2}u_{2}%
)+O(\varepsilon^{3}).
\end{align*}
We define the nilpotent approximation as:%
\begin{align*}
(\mathcal{\hat{P}}_{5,2})\text{ \ \ }\dot{x}_{1}  & =u_{1},\dot{x}_{2}%
=u_{2},\dot{y}=(\frac{x_{2}}{2}u_{1}-\frac{x_{1}}{2}u_{2}),\\
\dot{z}  & =x_{2}(\frac{x_{2}}{2}u_{1}-\frac{x_{1}}{2}u_{2}),\\
\dot{w}  & =\delta(w)x_{1}(\frac{x_{2}}{2}u_{1}-\frac{x_{1}}{2}u_{2}%
).\text{\ }%
\end{align*}

The estimates necessary to reduce to Nilpotent approximation are:%
\begin{align}
||x(t)-\hat{x}(t)||  & \leq c\varepsilon^{4},||y(t)-\hat{y}(t)||\leq
c\varepsilon^{3},\label{ff2}\\
||z(t)-\hat{z}(t)||  & \leq c\varepsilon^{4},||w(t)-\hat{w}(t)||\leq
c\varepsilon^{4}.\nonumber
\end{align}

\textbf{Generic }$6-2$\textbf{\ case (proven in Appendix)}$:$%
\begin{align}
\dot{x}_{1}  & =u_{1}+0(\varepsilon^{3}),\dot{x}_{2}=u_{2}+0(\varepsilon
^{3}),\label{nf62}\\
\dot{y}  & =(\frac{x_{2}}{2}u_{1}-\frac{x_{1}}{2}u_{2})+O(\varepsilon
^{2}),\nonumber\\
\dot{z}_{1}  & =x_{2}(\frac{x_{2}}{2}u_{1}-\frac{x_{1}}{2}u_{2})+O(\varepsilon
^{3}),\nonumber\\
\dot{z}_{2}  & =x_{1}(\frac{x_{2}}{2}u_{1}-\frac{x_{1}}{2}u_{2})+O(\varepsilon
^{3}),\nonumber\\
\dot{w}  & =Q_{w}(x_{1},x_{2})(\frac{x_{2}}{2}u_{1}-\frac{x_{1}}{2}%
u_{2})+O(\varepsilon^{4}),\nonumber
\end{align}

where $Q_{w}(x_{1},x_{2})$ is a quadratic form in $x$ depending smoothly on
$w.$

\bigskip

We define the nilpotent approximation as:%
\begin{align}
(\mathcal{\hat{P}}_{6,2})\text{ \ \ }\dot{x}_{1}  & =u_{1},\dot{x}_{2}%
=u_{2},\dot{y}=(\frac{x_{2}}{2}u_{1}-\frac{x_{1}}{2}u_{2}),\label{nil62}\\
\dot{z}_{1}  & =x_{2}(\frac{x_{2}}{2}u_{1}-\frac{x_{1}}{2}u_{2}),\dot{z}%
_{2}=x_{1}(\frac{x_{2}}{2}u_{1}-\frac{x_{1}}{2}u_{2}),\nonumber\\
\text{\ }\dot{w}  & =Q_{w}(x_{1},x_{2})(\frac{x_{2}}{2}u_{1}-\frac{x_{1}}%
{2}u_{2}).\nonumber
\end{align}

The estimates necessary to reduce to Nilpotent approximation are:%
\begin{align}
||x(t)-\hat{x}(t)||  & \leq c\varepsilon^{4},||y(t)-\hat{y}(t)||\leq
c\varepsilon^{3},\label{fff3}\\
||z(t)-\hat{z}(t)||  & \leq c\varepsilon^{4},||w(t)-\hat{w}(t)||\leq
c\varepsilon^{5}.\nonumber
\end{align}

In fact, the proof given in Appendix, of the reduction to this normal form,
contains the other cases 4-2 and 5-2.

\subsubsection{\bigskip Invariants in the 6-2 case, and the ball with a
trailer}

Let us consider a one form $\omega$ that vanishes on $\Delta^{\prime\prime
}=[\Delta,[\Delta,\Delta]].\ $Set $\alpha=d\omega_{|\Delta},$ the restriction
of $d\omega$ to $\Delta.$ Set $H=[F_{1},F_{2}],$ $I=[F1,H], $ $J=[F_{2},H],$
and consider the $2\times2$ matrix $A(\xi)=\left(
\begin{array}
[c]{cc}%
d\omega(F_{1},I) & d\omega(F_{2},I)\\
d\omega(F_{1},J) & d\omega(F_{2},J)
\end{array}
\right)  .$

Due tu Jacobi Identity,\ $\ A(\xi)$ is a symmetric matrix.\ It is also equal
to $\left(
\begin{array}
[c]{cc}%
\omega([F_{1},I]) & \omega([F_{2},I])\\
\omega([F_{1},J]) & \omega([F_{2},J])
\end{array}
\right)  ,$ using the fact that $\omega([X,Y])=d\omega(X,Y)$ in restriction to
$\Delta^{\prime\prime}.$

Let us consider a gauge transformation, i.e. a feedback that preserves the
metric (i.e. a change of othonormal frame $(F_{1},F_{2})$ obtained by setting
$\tilde{F}_{1}=\cos(\theta(\xi))F_{1}+\sin(\theta(\xi))F_{2\text{ }},$
$\tilde{F}_{2}=-\sin(\theta(\xi))F_{1}+\cos(\theta(\xi))F_{2\text{ }}).$

It is just a matter of tedious computations to check that the matrix $A(\xi) $
is changed for $\tilde{A}(\xi)=R_{\theta}A(\xi)R_{-\theta}.$ On the other
hand, the form $\omega$ is defined modulo muttiplication by a nonzero function
$f(\xi),$ and the same holds for $\alpha,$ since $d(f\omega)=fd\omega
+df\wedge\omega,$ and $\omega$ vanishes over $\Delta^{\prime\prime}.$ The
following lemma follows:

\begin{lemma}
\label{lem62inv}The ratio $r(\xi)$ of the (real) eigenvalues of $A(\xi)$ is an
invariant of the structure.
\end{lemma}

Let us now consider the normal form (\ref{nf62}), and compute the form
$\omega=\omega_{1}dx_{1}+...+\omega_{6}dw$ along $\Gamma$ (that is, where
$x,y,z=0).$ Computing all the brackets show that $\omega_{1}=\omega
_{2}=...=\omega_{5}=0.$ This shows also that in fact, along $\Gamma$, $A(\xi)$
is just the matrix of the quadratic form $Q_{w}.\ $We get the following:

\begin{lemma}
\label{lem62inv1} The invariant $r(\Gamma(t))$ of the problem $\mathcal{P}$ is
the same as the invariant $\hat{r}(\Gamma(t))$ of the nilpotent approximation
along $\Gamma.$
\end{lemma}

Let us compute the ratio $r$ for the ball with a trailer, \ Equation
(\ref{balltrailer}). We denote by $A_{1},A_{2}$ the two right-invariant vector
fields over $So(3,\mathbb{R)}$ appearing in (\ref{balltrailer}).\ We have:%

\begin{align*}
F_{1}  & =\frac{\partial}{\partial x_{1}}+A_{1}-\frac{1}{L}\cos(\theta
)\frac{\partial}{\partial\theta},\\
F_{2}  & =\frac{\partial}{\partial x_{2}}+A_{2}-\frac{1}{L}\sin(\theta
)\frac{\partial}{\partial\theta}.\\
\lbrack A_{1},A_{2}]  & =A_{3},[A_{1},A_{3}]=-A_{2},[A_{2},A_{3}]=A_{1}.
\end{align*}

Then, we compute the brackets: $H=A_{3-}\frac{1}{L^{2}}\frac{\partial
}{\partial\theta},$ $I=-A_{2-}\frac{1}{L^{3}}\sin(\theta)\frac{\partial
}{\partial\theta},$ $J=A_{1+}\frac{1}{L^{3}}\cos(\theta)\frac{\partial
}{\partial\theta},$ $[F_{1},I]=-A_{3-}\frac{1}{L^{4}}\frac{\partial}%
{\partial\theta},$ $[F_{1},J]=0=[F_{2},I],$ $[F_{2},J]=-A_{3-}\frac{1}{L^{4}%
}\frac{\partial}{\partial\theta}.$ Then:

\begin{lemma}
\label{balltrailerratio}For the ball with a trailer, the ratio $r(\xi)=1.$
\end{lemma}

These two last lemmas are a key point in the section \ref{new}: theyl imply in
particular that the system of geodesics of the nilpotent approximation is
integrable in Liouville sense, as we shall see.

\section{Results\label{results}}

In this section, we summarize and comment most of the results obtained in the
papers \cite{GM,GZ,GZ2,GZ3,ano,jak}.

\subsection{General results\label{genres}}

We need the concept of an $\varepsilon$-modification of an asymptotic optimal synthesis.

\begin{definition}
Given a one parameter family of (absolutely continuous, arclength
parametrized) admissible curves $\gamma_{\varepsilon}:$ $[0,T_{\gamma
_{\varepsilon}}]\rightarrow\mathbb{R}^{n},$ \textbf{an }$\varepsilon
$\textbf{-modification of }$\gamma_{\varepsilon}$ is another one parameter
family of (absolutely continuous, arclength parametrized) admissible curves
$\tilde{\gamma}_{\varepsilon}:$ $[0,T_{\tilde{\gamma}_{\varepsilon}%
}]\rightarrow\mathbb{R}^{n}$ such that for all $\varepsilon$ and for some
$\alpha>0$, if $[0,T_{\gamma_{\varepsilon}}]$ is splitted into subintervals of
length $\varepsilon,$ $[0,\varepsilon],$ $[\varepsilon,2\varepsilon],$
$[2\varepsilon,3\varepsilon],...$ then:

1. $[0,T_{\tilde{\gamma}_{\varepsilon}}]$ is splitted into corresponding
intervals, $[0,\varepsilon_{1}],$ $[\varepsilon_{1},\varepsilon_{1}%
+\varepsilon_{2}],$ $[\varepsilon_{1}+\varepsilon_{2},\varepsilon
_{1}+\varepsilon_{2}+\varepsilon_{3}],...$ with $\varepsilon\leq$
$\varepsilon_{i}<\varepsilon(1+\varepsilon^{\alpha}),$ $i=1,2,...,$

2.\ for each couple of an interval $I_{1}=[\tilde{\varepsilon}_{i}%
,\tilde{\varepsilon}_{i}+\varepsilon],$ (with $\tilde{\varepsilon}_{0}=0,$
$\tilde{\varepsilon}_{1}=\varepsilon_{1},$ $\tilde{\varepsilon}_{2}%
=\varepsilon_{1}+\varepsilon_{2},...$) and the respective interval
$I_{2}=[i\varepsilon,(i+1)\varepsilon],$ $\frac{d}{dt}(\tilde{\gamma})$ and
$\frac{d}{dt}(\gamma)$ coincide over $I_{2},$ i.e.:%
\[
\frac{d}{dt}(\tilde{\gamma})(\tilde{\varepsilon}_{i}+t)=\frac{d}{dt}%
(\gamma)(i\varepsilon+t),\text{ for almost all }t\in\lbrack0,\varepsilon].
\]

\end{definition}

\begin{remark}
This concept of an $\varepsilon$\textbf{-modification}\ is for the following
use: we will construct asymptotic optimal syntheses for the nilpotent
approximation $\mathcal{\hat{P}}$ of problem $\mathcal{P}$. Then, the
asymptotic optimal syntheses have to be slightly modified in order to realize
the interpolation constraints for the original (non-modified) problem. This
has to be done "slightly" for the length of paths remaining equivalent.

In this section it is always assumed but not stated that \textbf{we consider
generic problems only}.\ One first result is the following:
\end{remark}

\begin{theorem}
\label{eqnil}In the cases 2-step bracket generating, 4-2, 5-2, 6-2, (without
singularities), an asymptotic optimal synthesis [relative to the entropy] for
$\mathcal{P}$ is obtained as an $\varepsilon$-modification of an asymptotic
optimal synthesis for the nilpotent approximation $\mathcal{\hat{P}}.$ As a
consequence the entropy $E(\varepsilon)$ of $\mathcal{P}$ is equal to the
entropy $\hat{E}(\varepsilon)$ of $\mathcal{\hat{P}}.$
\end{theorem}

This theorem is proven in \cite{GZ3}.\ However, we can easily get an idea of
the proof, using the estimates of formulas (\ref{ff0}, \ref{ff1}, \ref{ff2},
\ref{fff3}).

All these estimates show that, if we apply an $\varepsilon$-interpolating
strategy to $\mathcal{\hat{P}}$, and the same controls to $\mathcal{P}$, at
time $\varepsilon$ (or length $\varepsilon$-since it is always possible to
consider arclength-parametrized trajectories), the enpoints of the two
trajectories are at subriemannian distance (either $d$ or $\hat{d})$ of order
$\varepsilon^{1+\alpha},$ for some $\alpha>0.\ $Then the contribution to the
entropy of $\mathcal{P}$, due to the correction necessary to interpolate
$\Gamma$ will have higher order.

\bigskip

Also, in the one-step bracket-generating case, we have the following equality:

\begin{theorem}
\label{2pi}(one step bracket-generating case, corank $k\leq3)$ The entropy is
equal to $2\pi$ times the metric complexity: $E(\varepsilon)=2\pi
MC(\varepsilon).$
\end{theorem}

The reason for this distinction between corank less or more than 3 is very
important, and will be explained in the section \ref{onestep}.

Another very important result is the following \textbf{logarithmic lemma},
that describes what happens in the case of a (generic) singularity of
$\Delta.\ $\ In the absence of such singularities, as we shall see, we shall
always have formulas of the following type, for the entropy (the same for the
metric complexity):
\begin{equation}
E(\varepsilon)\simeq\frac{1}{\varepsilon^{p}}%
{\displaystyle\int\limits_{\Gamma}}
\frac{dt}{\chi(t)},\label{ent1}%
\end{equation}

\bigskip

where $\chi(t)$ is a certain invariant along $\Gamma$.\ When the curve
$\Gamma(t)$ crosses tranversally a codimension-1 singularity (of
$\Delta^{\prime},$ or $\Delta^{\prime\prime}),$ the invariant $\chi(t)$
vanishes.\ This may happen at isolated points $t_{i},$ \ $i=1,...r.$ In that
case, we always have the following:

\begin{theorem}
\label{logl}(logarithmic lemma).\ The entropy (resp. the metric complexity)
satisfies:%
\[
E(\varepsilon)\simeq-2\frac{\ln(\varepsilon)}{\varepsilon^{p}}\sum_{i=1}%
^{r}\frac{1}{\rho(t_{i})},\text{ \ \ where }\rho(t)=|\frac{d\chi(t)}{dt}|.
\]

\end{theorem}

On the contrary, there are also generic codimension 1 singularities where the
curve $\Gamma$, at isolated points, becomes tangent to $\Delta,$ or
$\Delta^{\prime},...$ At these isolated points, the invariant $\chi(t)$ of
Formula \ref{ent1} tends to infinity.\ In that case, \textbf{the formula
\ref{ent1} remains valid} (the integral converges).

\subsection{Generic distribution in $\mathbb{R}^{3}$\label{contact3}}

This is the simplest case, and is is important, since many cases just reduce
to it.\ Let us describe it in details.

Generically, the 3-dimensional space $M$ contains a 2-dimensional singularity
(called the Martinet surface, denoted by $\mathcal{M)}.$ This singularity is a
smooth surface, and (except at isolated points on $\mathcal{M)},$ the
distribution $\Delta$ is not tangent to $\mathcal{M}.$ Generically, the curve
$\Gamma$ crosses $\mathcal{M}$ transversally at a finite number of isolated
points $t_{i},$ $i=1,...,r,.\ $These points are not the special isolated
points where $\Delta$ is tangent to $\mathcal{M}$ (this would be not generic).
They are called Martinet points. This number $r$ can be zero. Also, there are
other isolated points $\tau_{j},$ $j=1,...,l,$ at which $\Gamma$ is tangent to
$\Delta$ (which means that $\Gamma$ is almost admissible in a neighborhood of
$\tau_{j}).$ Out of $\mathcal{M}$, the distribution $\Delta$ is a contact
distribution (a generic property).

Let $\omega$ be a one-form that vanishes on $\Delta$ and that is 1 on
$\dot{\Gamma}$, defined up to multiplication by a function which is 1 along
$\Gamma.$ Along $\Gamma,$ the restriction 2-form $d\omega_{|\Delta}$ can be
made into a skew-symmetric endomorphism $A(\Gamma(t))$ of $\Delta$ (skew
symmetric with respect to the scalar product over $\Delta),$ by duality:
$<A(\Gamma(t))X,Y>=d\omega(X,Y).$ Let $\chi(t)$ denote the moduli of the
eigenvalues of $A(\Gamma(t)).$ We have the following:

\begin{theorem}
\label{dim3}1.\ If $r=0,$ $MC(\varepsilon)\simeq\frac{2}{\varepsilon^{2}}%
{\displaystyle\int\limits_{\Gamma}}
\frac{dt}{\varkappa(t)}.\ $At points where $\chi(t)\rightarrow+\infty,$ the
formula is convergent.

2.\ If $r\neq0,$ $MC(\varepsilon)\simeq-2\frac{\ln(\varepsilon)}%
{\varepsilon^{2}}\sum_{i=1}^{r}\frac{1}{\rho(t_{i})},$ \ \ where
$\rho(t)=|\frac{d\chi(t)}{dt}|.$

3.\ $E(\varepsilon)=2\pi MC(\varepsilon).$
\end{theorem}

\begin{figure}[ptb]%
\centering
\includegraphics[
natheight=2.354000in,
natwidth=3.604500in,
height=2.1075in,
width=3.2145in
]%
{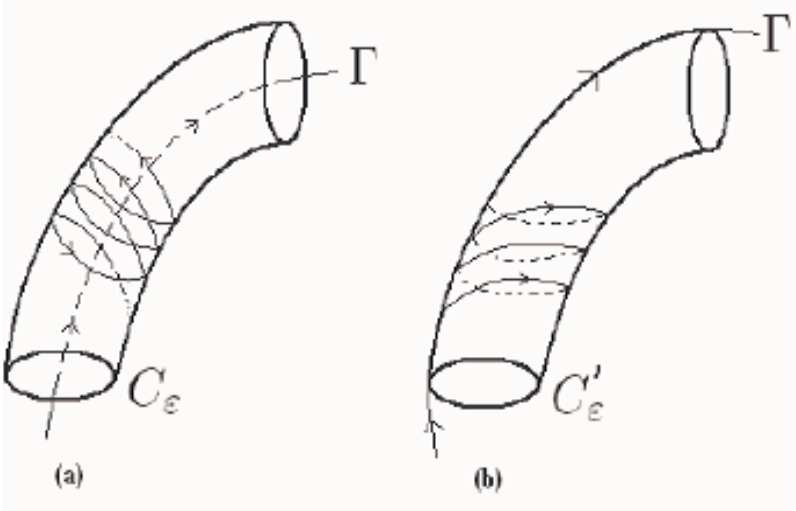}%
\caption{3-dimensional contact case}%
\label{Figcontact}%
\end{figure}

\begin{figure}[ptb]%
\centering
\includegraphics[
natheight=2.865100in,
natwidth=4.371600in,
height=2.0254in,
width=3.0787in
]%
{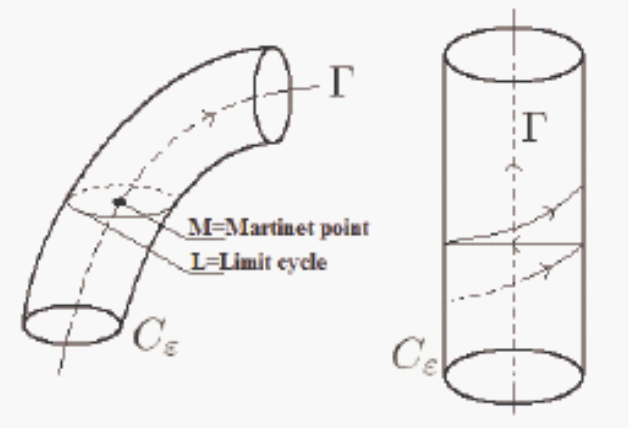}%
\caption{3-dimensional Martinet case}%
\label{figmar}%
\end{figure}

Let us describe the asymptotic optimal syntheses. They are shown on Figures
\ref{Figcontact}, \ref{figmar}.

\bigskip%

\bigskip

Figure \ref{Figcontact} concerns the case $r=0$ (everywhere contact type). The
points where the distribution $\Delta$ is not transversal to $\Gamma$ are
omitted (they again do not change anything)$.$ Hence $\Delta$ is also
transversal to the cylinders $C_{\varepsilon}$, for $\varepsilon$
small$.\ $Therefore, $\Delta$ defines (up to sign) a vector field
$X_{\varepsilon}$ on $C\varepsilon,$ tangent to $\Delta,$ that can be chosen
of length 1. The asymptotic optimal synthesis consists of: 1.\ Reaching
$C_{\varepsilon}$ from $\Gamma(0),$ 2.\ Follow a trajectory of $X_{\varepsilon
},$ 3.\ Join $\Gamma(t).\ $The steps 1 and 3 cost 2$\varepsilon,$ which is
neglectible w.r.t.\ the full metric complexity.\ To get the optimal synthesis
for the interpolation entropy, one has to make the same construction, but
starting from a subriemannian cylinder $C_{\varepsilon}^{\prime}$ tangent to
$\Gamma.$

In normal coordinates, in that case, the $x$-trajectories are just circles,
and the corresponding optimal controls are just trigonometric functions, with
period $\frac{2\pi}{\varepsilon}.$

Figure \ref{figmar} concerns the case $r\neq0$ (crossing Martinet
surface).\ At a Martinet point, the vector-field $X_{\varepsilon}$ has a limit
cycle, which is not tangent to the distribution. The asymptotic optimal
strategy consists of: a.\ following a trajectory of $X_{\varepsilon} $ till
reaching the height of the center of the limit cycle, b.\ crossing the
cylinder, with a neglectible cost $2\varepsilon,$ c.\ Following a trajectory
of the opposite vector field $-X_{\varepsilon}.$ The strategy for entropy is
similar, but using the tangent cylinder $C_{\varepsilon}^{\prime}.$%

\subsection{The one-step bracket-generating case\label{onestep}}

For the corank $k\leq3,$ the situation is very similar to the 3-dimensional
case.\ It can be competely reduced to it.\ For details, see \cite{GZ2}.

At this point, this strange fact appears: there is the limit corank $k=3.$\ If
$k>3$ only, new phenomena appear. Let us explain now the reason for this$\ $

Let us consider the following mapping $\mathcal{B}_{\xi}:\Delta_{\xi}%
\times\Delta_{\xi}\rightarrow T_{x}M/\Delta_{\xi},$ $(X,Y)\rightarrow\lbrack
X,Y]+\Delta_{\xi}.\ $It is a well defined \textbf{tensor mapping} , which
means that it actually applies to vectors (and not to vector fields, as
expected from the definition). This is due to the following formula, for a
one-form $\omega:$ $\ d\omega(X,Y)=\omega([X,Y])+\omega(Y)X-\omega(X)Y.$ Let
us call $I_{\xi}$ the image by $\mathcal{B}_{\xi}$ of the product of two unit
balls in $\Delta_{\xi}.\ $The following holds:

\begin{theorem}
\label{convexity} For a generic $\mathcal{P}$, for $k\leq3,$ the sets
$I_{\Gamma(t)}$ are \textbf{convex}.
\end{theorem}

This theorem is shown in \cite{GZ2}, with the consequences that we will state
just below.

This is no more true for $k>3,$ the first catastrophic case being the case
10-4 (a $p=4$ distribution in $\mathbb{R}^{10}).$ The intermadiate cases
$k=4,5$ in dimension 10 are interesting, since on some open subsets of
$\Gamma,$ the convexity property may hold or not. These cases are studied in
the paper \cite{ano}.

The main consequence of this convexity property is that everything reduces
(out of singularities where the logarithmic lemma applies) to the
3-dimensional contact case, as is shown in the paper \cite{GZ2}.\ We briefly
summarize the results.

Consider the one forms $\omega$ that vanish on $\Delta$ and that are 1 on
$\dot{\Gamma},$ and again, by the duality w.r.t. the metric over $\Delta,$
define $d\omega_{|\Delta}(X,Y)=<AX,Y>,$ for vector fields $X,Y$ in $\Delta.$
Now, we have \ along $\Gamma,$ a ($k-1)$-parameter affine family of skew
symmetric endomorphisms $A_{\Gamma(t)}$ of $\Delta_{\Gamma(t).}$ Say,
$A_{\Gamma(t)}(\lambda)=A_{\Gamma(t)}^{0}+%
{\displaystyle\sum\limits_{i=1}^{k-1}}
\lambda_{i}A_{\Gamma(t)}^{i}.$ Set $\chi(t)=\inf_{\lambda}||A_{\Gamma
(t)}(\lambda)||=||A_{\Gamma(t)}(\lambda^{\ast}(t))||.$

Out of isolated points of $\Gamma$ (that count for nothing in the metric
complexity or in the entropy), the $t-$one parameter family $A_{\Gamma
(t)}(\lambda^{\ast}(t))$ can be smoothly block-diagonalized (with $2\times2$
bloks), using a gauge transformation along $\Gamma$. After this gauge
transformation, the 2-dimensional eigenspace corresponding to the largest (in
moduli) eigenvalue of $A_{\Gamma(t)}(\lambda^{\ast}(t)),$ corresponds to the
two first coordinates in the distribution, and to the 2 first controls.\ In
the asymptotic optimal synthesis, all other controls are put to zero [here the
convexity property is used], and the picture of the asymptotic optimal
synthesis is exactly that of the 3-dimensional contact case. We still have the formulas:%

\[
MC(\varepsilon)\simeq\frac{2}{\varepsilon^{2}}%
{\displaystyle\int\limits_{\Gamma}}
\frac{dt}{\varkappa(t)},\text{ \ \ }E(\varepsilon)=2\pi MC(\varepsilon).
\]

The case $k>3$ was first treated in \cite{GZ3} in the 10-dimensional case, and
was completed in general in \cite{jak}.

\bigskip

In that case, the situation does not reduce to the 3-dimensional contact case:
the optimal controls, in the asymptotic optimal synthesis for the nilpotent
approximation are still trigonometric controls, but with different periods
that are successive integer multiples of a given basic period. New invariants
$\lambda_{\theta(t)}^{j}$ appear, and the formula for the entropy is:%
\[
E(\varepsilon)\simeq\frac{2\pi}{\varepsilon^{2}}\int_{0}^{T}\frac{\sum
_{j=1}^{r}j\lambda_{\theta}^{j}}{\sum_{j=1}^{r}(\lambda_{\theta}^{j})^{2}%
}d\theta,
\]
the optimal controls being of the form:%

\begin{align}
u_{2j-1}(t) &  =-\sqrt{\frac{j\lambda_{\theta(t)}^{j}}{\sum_{j=1}^{r}%
j\lambda_{\theta(t)}^{j}}}\sin(\frac{2\pi jt}{\varepsilon}),\label{controls}\\
u_{2j}(t) &  =\sqrt{\frac{j\lambda_{\theta(t)}^{j}}{\sum_{j=1}^{r}%
j\lambda_{\theta(t)}^{j}}}\cos(\frac{2\pi jt}{\varepsilon}),\text{
\ \ }j=1,...,r\nonumber\\
{u_{2r+1}(t)} &  =0{\text{ if }p\text{ is odd }}.\nonumber
\end{align}

These last formulas hold in the free case only (i.e. the case where the corank
$k=\frac{p(p-1)}{2},$ the dimension of he second homogeneous component of the
free Lie-algebra with $p$ generators).\ The non free case is more complicated
(see \cite{jak}).

To prove all the results in this section, one has to proceed as follows: 1.
use the theorem of reduction to nilpotent approximation (\ref{eqnil}), and
2.\ use the Pontriaguin'smaximum principle on the normal form of the nilpotent
approximation, in normal coordinates

\subsection{The 2-control case, in $\mathbb{R}^{4}$ and $\mathbb{R}^{5}.$}

These cases correspond respectively to the car with a trailer (Example
\ref{ctrl}) and the ball on a plate (Example \ref{bpln}).

We use also the theorem \ref{eqnil} of reduction to Nilpotent approximation,
and we consider the normal forms $\mathcal{\hat{P}}_{4,2},$ $\mathcal{\hat{P}%
}_{5,2}$ of Section \ref{nffs}. In both cases, we change the variable $w$ for
$\tilde{w}$ such that $d\tilde{w}=\frac{dw}{\delta(w)}.\ $\ We look for
arclength-parametrized trajectories of the nilpotent approximation (i.e.
$(u_{1})^{2}+(u_{2})^{2}=1),$ that start from $\Gamma(0),$ and reach $\Gamma$
in fixed time $\varepsilon,$ maximizing $%
{\displaystyle\int\limits_{0}^{\varepsilon}}
\dot{w}(\tau)d\tau.$ Abnormal extremals do no come in the picture, and optimal
curves correspond to the hamiltonian
\[
H=\sqrt{(PF_{1})^{2}+(PF_{2})^{2}},
\]
where $P$ is the adjoint vector.\ It turns out that, in our normal
coordinates, the same trajectories are optimal for both the $4$-$2$ and the
$5$-$2$ case (one has just to notice that the solution of the $4$-$2$ case
meets the extra interpolation condition corresponding to the 5-2 case).\ 

Setting as usual $u_{1}=\cos(\varphi)=PF_{1},u_{2}=\sin(\varphi)=PF_{2},$ we
get $\dot{\varphi}=P[F_{1},F_{2}],\ddot{\varphi}=-P[F_{1},[F_{1},F_{2}%
]]PF_{1}-P[F_{2},[F_{1},F_{2}]]PF_{2}.$

At this point, we have to notice that only the components $P_{x_{1}},P_{x_{2}%
}$ of the adloint vector $P$ are not constant (the hamiltonian in the
nilpotent approximation depends only on the $x$-variables), $\ $therefore,
$P[F_{1},[F_{1},F_{2}]]$ and $P[F_{2},[F_{1},F_{2}]]$ are constant (the third
brackets are also constant vector fields).\ Hence, $\ddot{\varphi}=\alpha
\cos(\varphi)+\beta\sin(\varphi)$ $=\alpha\dot{x}_{1}+\beta\dot{x}_{2}$ for
appropriate constants $\alpha,\beta.$ It follows that, for another constant
$k,$ we have, for the optimal curves of the nilpotent approximation, in normal
coordinates $x_{1},x_{2}:$%
\begin{align*}
\dot{x}_{1}  & =\cos(\varphi),\dot{x}_{2}=\sin(\varphi),\\
\dot{\varphi}  & =k+\lambda x_{1}+\mu x_{2}.
\end{align*}

\begin{remark}
\label{remcurv}1.\ It means that we are looking for curves in the $x_{1}%
,x_{2}$ plane, whose curvature is an affine function of the position,

2.\ In the two-step bracket generating case (contact case), otimal curves were
circles, i.e. curves of constant curvature,

3.\ the conditions of $\varepsilon$-interpolation of $\Gamma$ say that these
curves must be periodic (there will be more details on this point in the next
section), that the area of a loop must be zero $(y(\varepsilon)=0), $ and
finally (in the 5-2 case) that another moment must be zero.
\end{remark}

\begin{figure}[ptb]%
\centering
\includegraphics[
natheight=1.604200in,
natwidth=3.989400in,
height=1.3872in,
width=3.4143in
]%
{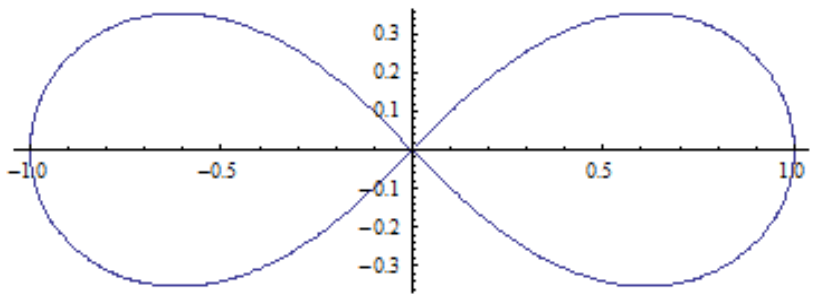}%
\caption{The dance of minimum entropy, for 3rd bracket.}%
\label{elastica}%
\end{figure}

It is easily seen that such a curve, meeting these interpolation conditions,
must be an elliptic curve of elastica-type. The periodicity and vanishing
surface requirements imply that it is the only periodic elastic curve shown on
Figure \ref{elastica}, parametrized in a certain way.%

The formulas are, in terms of the standard Jacobi elliptic functions:%
\begin{align*}
u_{1}(t)  & =1-2dn(K(1+\frac{4t}{\varepsilon}))^{2},\\
u_{2}(t)  & =-2dn(K(1+\frac{4t}{\varepsilon}))sn(K(1+\frac{4t}{\varepsilon
}))\sin(\frac{\varphi_{0}}{2}),
\end{align*}
where $\varphi_{0=130%
{{}^\circ}%
}$ (following \cite{Love}, p.\ 403) and $\varphi_{0}=130,692%
{{}^\circ}%
$ following Mathematica$^{\textregistered},$ with $k=\sin(\frac{\varphi_{0}%
}{2})$ and $K(k)$ is the quarter period of the Jacobi elliptic functions. The
trajectory on the $x_{1},x_{2}$ plane, shown on Figure \ref{elastica}, has
equations:
\begin{align*}
x_{1}(t)  & =-\frac{\varepsilon}{4K}[\frac{-4Kt}{\varepsilon}+2(Eam(\frac
{4Kt}{\varepsilon}+K)-Eam(K))],\\
x_{2}(t)  & =k\frac{\varepsilon}{2K}cn(\frac{4Kt}{\varepsilon}+K).
\end{align*}

On the figure \ref{fig2}, one can clearly see, at the contact point of the
ball with the plane, a trajectory which is a "repeated small deformation" of
this basic trajectory.

The formula for the entropy is, in both the 4-2 and 5-2 cases:%
\[
E(\varepsilon)=\frac{3}{2\sigma\varepsilon^{3}}\int_{\Gamma}\frac{dt}%
{\delta(t)},
\]
where $\sigma$ is a universal constant, $\sigma\approx0.00580305.$

Details of computations on the 4-2 case can be found in \cite{GZ3}, and in
\cite{ano} for the 5-2 case.

\section{The ball with a trailer\label{new}}

We start by using Theorem \ref{eqnil}, to reduce to the nilpotent
approximation along $\Gamma:$%

\begin{align}
(\mathcal{\hat{P}}_{6,2})\text{ \ \ }\dot{x}_{1}  & =u_{1},\dot{x}_{2}%
=u_{2},\dot{y}=(\frac{x_{2}}{2}u_{1}-\frac{x_{1}}{2}u_{2}),\\
\dot{z}_{1}  & =x_{2}(\frac{x_{2}}{2}u_{1}-\frac{x_{1}}{2}u_{2}),\dot{z}%
_{2}=x_{1}(\frac{x_{2}}{2}u_{1}-\frac{x_{1}}{2}u_{2}),\nonumber\\
\text{\ }\dot{w}  & =Q_{w}(x_{1},x_{2})(\frac{x_{2}}{2}u_{1}-\frac{x_{1}}%
{2}u_{2}).\nonumber
\end{align}

By Lemma \ref{balltrailerratio}, we can consider that
\begin{equation}
Q_{w}(x_{1},x_{2})=\delta(w)((x_{1})^{2}+(x_{2})^{2})\label{mainf}%
\end{equation}
where $\delta(w)$ is \textbf{the main invariant}. In fact, it is the only
invariant for the nilpotent approximation along $\Gamma.$ Moreover, if we
reparametrize $\Gamma$ by setting $dw:=\frac{dw}{\delta(w)},$ we can consider
that $\delta(w)=1.$

Then, we want to maximaize $\int\dot{w}dt$ in fixed time $\varepsilon,$ with
the interpolation conditions: $x(0)=0,y(0)=0,z(0)=0,w(0)=0,$ $x(\varepsilon
)=0,y(\varepsilon)=0,z(\varepsilon)=0.$

\bigskip

From Lemma \ref{periodl} in the appendix, we know that the optimal trajectory
is smooth and periodic, (of period $\varepsilon).$

Clearly, the optimal trajectory has also to be a length minimizer, then we
have to consider the usual hamiltonian for length: $H=\frac{1}{2}%
((P.F_{1})^{2}+(P.F_{2})^{2}),$ in which $P=(p_{1},...,p_{6})$ is the adjoint
vector. It is easy to see that the abnormal extremals do not come in the
picture (cannot be optimal with our additional interpollation conditions), and
in fact, we will show that \textbf{the hamiltonian system corresponding to the
hamiltonian }$H$\textbf{\ is integrable}.

\begin{remark}
This integrability property is no more true in the general 6-2
case.\textbf{\ It holds only for the ball with a trailer.}
\end{remark}

As usual, we work in Poincar\'{e} coordinates, i.e. we consider level
$\frac{1}{2}$ of the hamiltonian $H,$ and we set:
\[
u_{1}=PF=\sin(\varphi),\text{ \ \ }u_{2}=PG=\cos(\varphi).
\]

Differentiating twice, we get
\[
\dot{\varphi}=P[F,G],\text{ }\ddot{\varphi}=-PFFG.PF-PGFG.PG,
\]
where $FFG=[F,[F,G]]$ and $GFG=[G,[F,G]].\ $\ We set $\lambda=-PFFG,$
$\mu=-PGFG.\ $We get that:%

\begin{equation}
\ddot{\varphi}=\lambda\sin(\varphi)+\mu\cos(\varphi).\label{phi2}%
\end{equation}

Now, we compute $\dot{\lambda}$ and $\dot{\mu}.$ We get, with similar
notations as above for the brackets (we bracket from the left):%

\begin{align*}
\dot{\lambda}  & =PFFFG.PF+PGFFG.PG,\\
\dot{\mu}  & =PFGFG.PF+PGGFG.PG,
\end{align*}
and computing the brackets, we see that $GFFG=FGFG=0.\ $Also, since the
hamiltonian does not depend on $y,z,w,$ we get that $p_{3},p_{4},p_{5},p_{6}$
are constants.\ Computing the brackets $FFG$ and $GFG$ , we get that
\[
\lambda=\frac{3}{2}p_{4}+p_{6}x_{1},\text{ \ }\mu=\frac{3}{2}p_{5}+p_{6}x_{2},
\]
and then, $\dot{\lambda}=p_{6}\sin(\varphi)$ and $\dot{\mu}=p_{6}\cos
(\varphi).$ Then, by (\ref{phi2}), $\ddot{\varphi}=\frac{\lambda\dot{\lambda}%
}{p_{6}}+\frac{\mu\dot{\mu}}{p_{6}},$ and finally:%
\begin{align}
\dot{x}_{1}  & =\sin(\varphi),\text{ \ \ }\dot{x}_{2}=\cos(\varphi
),\label{sysint}\\
\dot{\varphi}  & =K+\frac{1}{2p_{6}}(\lambda^{2}+\mu^{2}),\nonumber\\
\dot{\lambda}  & =p_{6}\sin(\varphi),\text{ \ }\dot{\mu}=p_{6}\cos
(\varphi).\nonumber
\end{align}

Setting $\omega=\frac{\lambda}{p_{6}},\delta=\frac{\mu}{p_{6}},$ we obtain:%
\begin{align*}
\dot{\omega}  & =\sin(\varphi),\text{ \ }\dot{\delta}=\cos(\varphi),\\
\dot{\varphi}  & =K+\frac{p_{6}}{2}(\omega^{2}+\delta^{2}).
\end{align*}

It means that the plane curve $(\omega(t),\delta(t))$ has a curvature which is
a quadratic function of the distance to the origin.\ Then, the optimal curve
($x_{1}(t),x_{2}(t))$ projected to the horizontal plane of the normal
coordinates has a curvature which is a quadratic function of the distance to
some point. Following the lemma (\ref{eqcur}) in the appendix, this system of
equations is integrable.

\bigskip Summarizing all the results, we get the following theorem.

\begin{theorem}
\label{mainth}(\textbf{asymptotic optimal synthesis for the ball with a
trailer}) The asymptotic optimal synthesis is an $\varepsilon$-modification of
the one of the nilpotent approximation, which has the following properties, in
projection to the horizontal plane $(x_{1},x_{2})$ in normal coordinates:

1.\ It is a closed smooth periodic curve, whose curvature is a quadratic
function of the position, and a function of the square distance to some point,

2.\ The area and the 2$^{nd}$ order moments $\int_{\Gamma}x_{1}(x_{2}%
dx_{1}-x_{1}dx_{2})$ and $\int_{\Gamma}x_{2}(x_{2}dx_{1}-x_{1}dx_{2})$ are zero.

3.\ The entropy is given by the formula: $E(\varepsilon)=\frac{\sigma
}{\varepsilon^{4}}\int_{\Gamma}\frac{dw}{\delta(w)},$ where $\delta(w)$ is the
main invariant from (\ref{mainf}), and $\sigma$ is a universal constant.
\end{theorem}

In fact we can go a little bit further to integrate explicitely the system
(\ref{sysint}).\ Set $\bar{\lambda}=\cos(\varphi)\lambda-\sin(\varphi)\mu,$
$\bar{\mu}=\sin(\varphi)\lambda+\cos(\varphi)\mu.\ $we get:%
\begin{align*}
\frac{d\bar{\lambda}}{dt}  & =-\bar{\mu}(K+\frac{1}{2p_{6}}(\bar{\lambda}%
^{2}+\bar{\mu}^{2})),\\
\frac{d\bar{\mu}}{dt}  & =p_{6}+\bar{\lambda}(K+\frac{1}{2p_{6}}(\bar{\lambda
}^{2}+\bar{\mu}^{2})).
\end{align*}

\begin{figure}[ptb]%
\centering
\includegraphics[
natheight=6.718700in,
natwidth=2.817600in,
height=2.7536in,
width=1.1623in
]%
{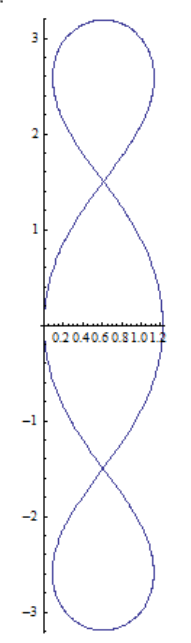}%
\caption{The dance of minimum entropy for the ball with a trailer}%
\label{fig62}%
\end{figure}

This is a 2 dimensional (integrable) hamiltonian system.\ The hamiltonian is:%
\[
H_{1}=-p_{6}\bar{\lambda}-\frac{2p_{6}}{4}(K+\frac{1}{2p_{6}}(\bar{\lambda
}^{2}+\bar{\mu}^{2}))^{2}.
\]
This hamiltonian system is therefore integrable, and solutions can be
expressed in terms of hyperelliptic functions. A liitle numerics now allows to
show, on figure \ref{fig62}, the optimal $x$-trajectory in the horizontal
plane of the normal coordinates.%

\begin{figure}[ptb]%
\centering
\includegraphics[
natheight=3.239600in,
natwidth=4.187400in,
height=2.4206in,
width=3.1228in
]%
{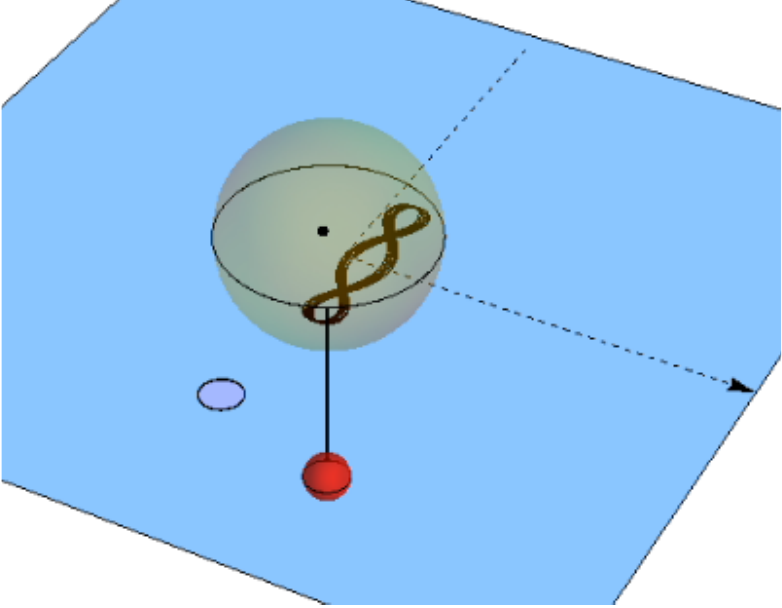}%
\caption{Parking the ball with a trailer}%
\label{figmov}%
\end{figure}

On the figure \ref{figmov}, we show the motion of the ball with a trailer on
the plane (motion of the contact point between the ball and the plane).Here,
the problem is to move along the $x$-axis, keeping constant the frame attached
to the ball and the angle of the trailer.%

\section{Expectations and conclusions\label{concl}}

Some movies of minimum entropy for the ball rolling on a plane and the ball
with a trailer are visible on the website ***************************.

\subsection{Universality of some pictures in normal coordinates}

Our first conclusion is the following: there are certain universal pictures
for the motion planning problem, in corank less or equal to 3, and in rank 2,
with 4 brackets at most (could be 5 brackets at a singularity, with the
logarithmic lemma).

These figures are, in the two-step bracket generating case: a circle, for the
third bracket, the periodic elastica, for the 4$^{th}$ bracket, the plane
curve of the figure \ref{fig62}.

They are periodic plane curves whose curvature is respectively: a constant, a
linear function of of the position, a quadratic function of the position.

\bigskip%
\begin{figure}[ptb]%
\centering
\includegraphics[
natheight=4.760800in,
natwidth=5.219100in,
height=2.9032in,
width=3.1816in
]%
{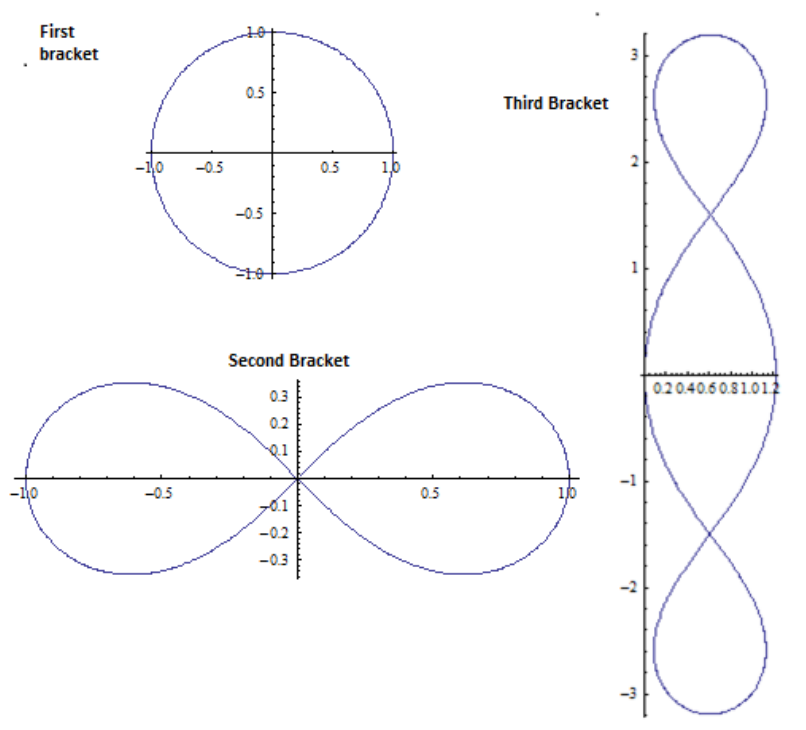}%
\caption{The universal movements in normal coordinates}%
\label{global}%
\end{figure}

This is, as shown on Figure \ref{global}, the clear beginning of a series.\ 

\subsection{Robustness}

As one can see, in many cases (2 controls, or corank $k\leq3),$ our strategy
is extremely robust in the following sense: the asymptotic optimal syntheses
do not depend, from the qualitative point of view, of the metric chosen.\ They
depend only on the number of brackets needed to generate the space.

\subsection{The practical importance of normal coordinates}

The main practical problem of implementation of our strategy comes with the
$\varepsilon$-modifications.\ How to compute them, how to implement?\ In fact,
the $\varepsilon$-modifications count at higher order in the entropy. But, if
not applied, they may cause deviations that are not neglectible.\ The high
order w.r.t. $\varepsilon$ in the estimates of the error between the original
system and its nilpotent approximation (Formulas \ref{ff0}, \ref{ff1},
\ref{ff2}, \ref{ff3}) make these deviations very small.\ It is why the use of
our concept of a nilpotent approximation along $\Gamma,$ based upon normal
coordinates is very efficient in practice.\ 

On the other hand, when a correction appears to be needed (after a
noneglectible deviation), it corresponds to brackets of lower order.\ For
example, in the case of the ball with a trailer (4$^{th}$ bracket), the
$\varepsilon$-modification corresponds to brackets of order 2 or 3.\ The
optimal pictures corresponding to these orders can still be used to perform
the $\varepsilon$-modifications.

\subsection{Final conclusion}

This approach, to approximate optimally nonadmissible paths of nonholomic
systems, looks very efficient, and in a sense, universal.\ Of course, the
theory is not complete, but the cases under consideration (first, 2-step
bracket-generating, and second, two controls) correspond to many practical
situations. But there is still a lot of work to do to in order to cover all
interesting cases.\ However, the methodology to go ahead is rather clear.

\section{\bigskip Appendix\label{app}}

\subsection{Appendix 1: Normal form in the 6-2 case\label{apnf62}}

We start from the general normal form (\ref{nf2}) in normal coordinates:%
\begin{align*}
\dot{x}_{1}  & =(1+(x_{2}^{{}})^{2}\beta)u_{1}-x_{1}x_{2}\beta u_{2},\text{
\ }\\
\text{\ }\dot{x}_{2}  & =(1+(x_{1}^{{}})^{2}\beta)u_{2}-x_{1}x_{2}\beta
u_{1},\\
\text{ \ }\dot{y}_{i}  & =(\frac{x_{2}}{2}u_{1}-\frac{x_{1}}{2}u_{2}%
)\gamma_{i}(y,w),\text{ }\\
\text{\ }\dot{w}  & =(\frac{x_{2}}{2}u_{1}-\frac{x_{1}}{2}u_{2})\delta(y,w)
\end{align*}

We will make a succession of changes of parametriztion of the surrface
$\mathcal{S}$ (w.r.t. which normal coordinates were constructed).\ These
coordinate changes will always preserve tha fact that $\Gamma(t)$ is the point
$x=0,y=0,w=t.$

Remind that $\beta$ vanishes on $\mathcal{S},$ and since $x$ has order $1,$ we
can already write on $T_{\varepsilon}$: $\dot{x}=u+O(\varepsilon^{3}).$
$\ $One of the $\gamma_{i}$'s (say $\gamma_{1})$ has to be nonzero (if not,
$\Gamma$ is tangent to $\Delta^{\prime}).$ Then, $y_{1}$ has order 2 on
$T_{\varepsilon}.$Set for $i>1,$ $\tilde{y}_{i}=y_{i}-\frac{\gamma_{i}}%
{\gamma_{1}}.\ $Differentiating, we get that $\frac{d\tilde{y}_{i}}{dt}%
=\dot{y}_{i}-\frac{\gamma_{i}}{\gamma_{1}}\dot{y}_{1}+O(\varepsilon^{2}),$ and
$z_{1}=\tilde{y}_{2},$ $z_{2}=\tilde{y}_{3}$ have order 3.\ We set also
$w:=w-\frac{\delta}{\gamma_{1}},$ and we are at the following point:%
\begin{align*}
\dot{x}  & =u+O(\varepsilon^{3}),\text{ \ }\dot{y}=(\frac{x_{2}}{2}u_{1}%
-\frac{x_{1}}{2}u_{2})\gamma_{1}(w)+O(\varepsilon^{2}),\\
\dot{z}_{i}  & =(\frac{x_{2}}{2}u_{1}-\frac{x_{1}}{2}u_{2})L_{i}%
(w).x+O(\varepsilon^{3}),\\
\dot{w}  & =(\frac{x_{2}}{2}u_{1}-\frac{x_{1}}{2}u_{2})\delta
(w).x+O(\varepsilon^{3}),
\end{align*}
where $L_{i}(w).x,$ $\delta(w).x$ are liner in $x.$ The function $\gamma
_{1}(w)$ can be put to 1 in the same way by setting $y:=\frac{y}{\gamma
_{1}(w)}.$ Now let $T(w)$ be an invertible 2$\times2$ matrix.\ Set $\tilde
{z}=T(w)z.$ It is easy to see that we can chose $T(w)$ for we get:%

\begin{align*}
\dot{x}  & =u+O(\varepsilon^{3}),\text{ \ }\dot{y}=(\frac{x_{2}}{2}u_{1}%
-\frac{x_{1}}{2}u_{2})+O(\varepsilon^{2}),\\
\dot{z}_{i}  & =(\frac{x_{2}}{2}u_{1}-\frac{x_{1}}{2}u_{2})x_{i}%
+O(\varepsilon^{3}),\\
\dot{w}  & =(\frac{x_{2}}{2}u_{1}-\frac{x_{1}}{2}u_{2})\delta
(w).x+O(\varepsilon^{3}),
\end{align*}

Another change of the form: $w:=w+L(w).x,$ where $L(w).x$ is linear in $x$
kills $\delta(w)$ and brings us to $\dot{w}=(\frac{x_{2}}{2}u_{1}-\frac{x_{1}%
}{2}u_{2})O(\varepsilon^{2}).\ $This $O(\varepsilon^{2})$ can be of the form
$Q_{w}(x)+h(w)y+O(\varepsilon^{3})$ where $Q_{w}(x)$ is quadratic in $x.$ If
we kill $h(w),$ we get the expected result.\ This is done with a change of
coordinates of the form: $w:=w+\varphi(w)\frac{y^{2}}{2}.$

\subsection{Appendix 2: Plane curves whose curvature is a function of the
distance to the origin\label{curvaturedist}}

This result was known already, see \cite{Singer}. However we provide here a
very simple proof.

Consider a plane curve $(x(t),y(t)),$ whose curvature is a function of the
distance from the origin, i.e.:%
\begin{equation}
\dot{x}=\cos(\varphi),\dot{y}=\sin(\varphi),\dot{\varphi}=k(x^{2}%
+y^{2}).\label{eqcur}%
\end{equation}

\bigskip Equation \ref{eqcur} is integrable.

\begin{proof}
Set $\bar{x}=x\cos(\varphi)+y\sin(\varphi),$ $\bar{y}=-x\sin(\varphi
)+y\cos(\varphi).$ Then $k(\bar{x}^{2}+\bar{y}^{2})=k(x^{2}+y^{2}).$ Just
computing, one gets:%

\begin{align}
\frac{d\bar{x}}{dt}  & =1+\bar{y}k(\bar{x}^{2}+\bar{y}^{2}),\label{ham}\\
\frac{d\bar{y}}{dt}  & =-\bar{x}k(\bar{x}^{2}+\bar{y}^{2}).\nonumber
\end{align}

We just show that (\ref{ham}) is a hamiltonian system.\ Since we are in
dimension 2, it is always Liouville-integrable. Then, we are looking for
solutions of the system of PDE's:%

\begin{align*}
\frac{\partial H}{\partial\bar{x}}  & =1+\bar{y}k(\bar{x}^{2}+\bar{y}^{2}),\\
\frac{\partial H}{\partial\bar{y}}  & =-\bar{x}k(\bar{x}^{2}+\bar{y}^{2}).
\end{align*}
But the Schwartz integrability conditions are satisfied: $\frac{\partial^{2}%
H}{\partial\bar{x}\partial\bar{y}}=\frac{\partial^{2}H}{\partial\bar
{y}\partial\bar{x}}=2\bar{x}\bar{y}k^{\prime}.$
\end{proof}

\subsection{Appendix 3: periodicity of the optimal curves in the 6-2 case
\label{periodicity}}

\begin{proof}
We consider the nilpotent approximation $\mathcal{\hat{P}}_{6,2}$ given in
formula \ref{nil62}:%

\begin{align}
(\mathcal{\hat{P}}_{6,2})\text{ \ \ }\dot{x}_{1}  & =u_{1},\dot{x}_{2}%
=u_{2},\dot{y}=(\frac{x_{2}}{2}u_{1}-\frac{x_{1}}{2}u_{2}),\\
\dot{z}_{1}  & =x_{2}(\frac{x_{2}}{2}u_{1}-\frac{x_{1}}{2}u_{2}),\dot{z}%
_{2}=x_{1}(\frac{x_{2}}{2}u_{1}-\frac{x_{1}}{2}u_{2}),\nonumber\\
\text{\ }\dot{w}  & =Q_{w}(x_{1},x_{2})(\frac{x_{2}}{2}u_{1}-\frac{x_{1}}%
{2}u_{2}).\nonumber
\end{align}

We consider the particular case of the ball with a trailer.\ Then, according
to Lemma \ref{balltrailerratio}, the ratio $r(\xi)=1.$

It follows that the last equation can be rewritten \ $\dot{w}=\delta
(w)((x_{1})^{2}+(x_{2})^{2})(\frac{x_{2}}{2}u_{1}-\frac{x_{1}}{2}u_{2})$ for
some never vanishing function $\delta(w)$ (vanishing would contradict the full
rank of $\Delta^{(4)}$). We can change the coordinate $w$ for $\tilde{w}$ such
that $d\tilde{w}=\frac{dw}{\delta(w)}.$

We get finally:%

\begin{align}
(\mathcal{\hat{P}}_{6,2})\text{ \ \ }\dot{x}_{1}  & =u_{1},\dot{x}_{2}%
=u_{2},\dot{y}=(\frac{x_{2}}{2}u_{1}-\frac{x_{1}}{2}u_{2}),\label{finalnil}\\
\dot{z}_{1}  & =x_{2}(\frac{x_{2}}{2}u_{1}-\frac{x_{1}}{2}u_{2}),\dot{z}%
_{2}=x_{1}(\frac{x_{2}}{2}u_{1}-\frac{x_{1}}{2}u_{2}),\nonumber\\
\ \dot{w}  & =((x_{1})^{2}+(x_{2})^{2})(\frac{x_{2}}{2}u_{1}-\frac{x_{1}}%
{2}u_{2})\nonumber
\end{align}

This is a right invariant system on $\mathbb{R}^{6}$ with cooordinates
$\xi=(\varsigma,w)=(x,y,z,w),$ for a certain Nilpotent Lie group structure
over $\mathbb{R}^{6}$ (denoted by $G).$ It is easily seen (just expressing
right invariance) that the group law is ot the form $(\varsigma_{2}%
,w_{2})(\varsigma_{1},w_{1})=$ $(\varsigma_{1}\ast\varsigma_{2},w_{1}%
+w_{2}+\Phi(\varsigma_{1},\varsigma_{2})),$ where $\ast$ is the multiplication
of another Lie group structure on $\mathbb{R}^{5},$ with coordinates
$\varsigma$ (denoted by $G_{0}).$ In fact, $G$ is a central extension of
$\mathbb{R}$ by $G_{0}.$
\end{proof}

\begin{lemma}
\label{periodl}The trajectories of (\ref{finalnil}) that maximize $%
{\displaystyle\int}
\dot{w}dt$ in fixed time $\varepsilon,$ with interpolating conditions
$\varsigma(0)=\varsigma(\varepsilon)=0,$ have a periodic projection on
$\varsigma$ (i.e.\ $\varsigma(t)$ is smooth and periodic of period
$\varepsilon).$
\end{lemma}

\begin{remark}
1.\ Due to the invariance with respect to the $w$ coordinate of
(\ref{finalnil}), it is equivalent to consider the problem with the more
restrictive terminal conditions $\varsigma(0)=\varsigma(\varepsilon)=0,$
$w(0)=0,$

2.\ The scheme of this proof works also to show periodicity in the 4-2 and 5-2 cases.
\end{remark}

The idea for the proof was given to us by A.\ Agrachev.

\begin{proof}
Let $(\varsigma,w_{1}),(\varsigma,w_{2})$ be initial and terminal points of an
optimal solution of our problem. By right translation by $(\varsigma^{-1},0),$
this trajectory is mapped into another trajectory of the system, with initial
and terminal points $(0,w_{1}+\Phi(\varsigma,\varsigma^{-1})) $ and
$(0,w_{1}+\Phi(\varsigma,\varsigma^{-1})).\ $Hence, this trajectory has the
same value of the cost $%
{\displaystyle\int}
\dot{w}dt.$ We see that the optimal cost is in fact independant of the
$\varsigma$-coordinate of the initial and terminal condition.\ 

Therefore, the problem is the same as maximizing $%
{\displaystyle\int}
\dot{w}dt$ but with the (larger) endpoint condition $\varsigma(0)=\varsigma
(\varepsilon)$ (free). Now, we can apply the general transversality conditions
of Theorem 12.15\ page 188 of \cite{AS}.\ It says that the initial and
terminal covectors $(p_{\varsigma}^{1},p_{w}^{1})$ and $(p_{\varsigma}%
^{2},p_{w}^{2})$ are such that $p_{\varsigma}^{1}=p_{\varsigma}^{2}.\ $This is
enough to show periodicity.
\end{proof}

\end{document}